\documentclass[]{interact}

\usepackage{epstopdf}
\usepackage[caption=false]{subfig}

\usepackage[numbers,sort&compress]{natbib}
\bibpunct[, ]{[}{]}{,}{n}{,}{,}
\renewcommand\bibfont{\fontsize{10}{12}\selectfont}
\makeatletter
\def\NAT@def@citea{\def\@citea{\NAT@separator}}
\makeatother

\theoremstyle{plain}
\newtheorem{theorem}{Theorem}[section]
\newtheorem{lemma}[theorem]{Lemma}

\theoremstyle{definition}
\newtheorem{definition}[theorem]{Definition}
\newtheorem{example}[theorem]{Example}

\usepackage{xcolor}

\theoremstyle{remark}
\newtheorem{remark}[theorem]{Remark}

\newcommand \+{^{\dagger }}
\begin{document}


\title{ Reverse Order Law for Closed Range Operators in Hilbert Spaces}

\author{\name{ Athira Satheesh K.\textsuperscript{a}, K. Kamaraj\textsuperscript{b} and P. Sam Johnson\textsuperscript{c} }
\affil{ \textsuperscript{a,}\textsuperscript{c}Department of Mathematical and Computational Sciences,
	National Institute of Technology Karnataka (NITK), Surathkal, Mangaluru 575 025, India\\ Email: athirachandri@gmail.com, sam@nitk.edu.in \\
	\textsuperscript{b}Department of Mathematics, University College of Engineering Arni, Thatchur, Arni\\ 632 326,
	India\\ Email: krajkj@yahoo.com}
}

\maketitle

\begin{abstract}
We present more than 50 results including some range inclusion results to characterize reverse order law for Moore-Penrose inverse of closed range Hilbert space operators. We use basic properties   of Moore-Penrose inverse to prove the results. Some examples are also provided to illustrate failure cases to hold the reverse order law in  infinite dimensional settings. 
\end{abstract}


\begin{amscode}47A05, 15A09.\end{amscode}

\section{Introduction}

One of the fundamental research problems in the theory of generalized inverses of matrices is to establish reverse order laws for generalized inverses of matrix products. It was Ivan Fredholm who seems to have first mentioned  the concept of generalized inverse in 1903. He formulated a pseudoinverse for a linear integral operator
which is not invertible in the ordinary sense.  Hilbert, Schmidt, Bounitzky, Hurwitz and other mathematicians had studied
the generalized inverses of integral operators and differential operators before  Moore introduced the generalized inverse of matrices by algebraic methods in 1920 \cite{moore}.
Bjerhammar rediscovered Moore's inverse and also noted the relationship
of generalized inverses to solutions of linear systems in 1951 \cite{arne}.
In 1955, Penrose \cite{penrose} extended Bjerhammar's
results and showed that  Moore's inverse for a given matrix $A$ is the unique matrix
$X$ satisfying the four equations:
\begin{eqnarray*}
AXA=A;
XAX=X;	
(AX)^* = AX;
(XA)^*= XA .
\end{eqnarray*}

In honour of Moore and Penrose this unique inverse is now commonly called
the Moore-Penrose inverse and denoted by $A^\dagger.$ Meanwhile, generalized inverses were defined for operators by Tseng \cite{Tseng}, Murray and von Neumann \cite{murray}, Nashed \cite{Nashed} and others.  Beutler discussed  generalized inverse for both bounded and unbounded operators with closed and arbitrary ranges  \cite{Beutler1,Beutler2}. Throughout the years the Moore-Penrose
inverse was extensively studied. One of the primary reasons for the same is its usefulness in
solving systems of linear equations, which constitutes an important application in various fields.

It is well known that the reverse order law $(AB)^{-1}=B^{-1}A^{-1} $ does not hold good for various generalized inverses such as Moore-Penrose inverse, Drazin inverse etc.  Cline attempted to find a reasonable representation for Moore-Penrose inverse of the product of  matrices \cite{cline} and  Greville found some necessary and sufficient conditions for reverse order law to hold in matrix settings \cite{Greville}.  The reverse order law problem for   bounded linear operators on
Hilbert spaces was analyzed by  Bouldin \cite{Bouldin1973, Bouldin1982} and Izumino \cite{Izumino}.  The theory of generalized inverses on infinite-dimensional
Hilbert spaces can be found in \cite{Israel, Wangbook, Harte}.

Tian \cite{Tian} established a group of rank formulas related to Moore-Penrose inverses  $(AB)^\dagger$ and  $B^\dagger A^\dagger.$ Djordjevi\'{c} and   Din\v{c}i\'{c}  extended Tian's result in \cite{Tian} to infinite dimensional settings using the matrix form of a bounded linear operator  \cite{Djordjevic2010}. Here, we also present extensions of Tians's results to infinite dimensional Hilbert space, but our approach is entirely different from Djordjevi\'{c} and   Din\v{c}i\'{c}. We used only basic properties of Moore-Penrose inverse and we claim that our proofs are much  simpler than Djordjevi\'{c} and   Din\v{c}i\'{c}.   In 1963, Arghiriade gave a characterization result for reverse order law to hold in finite dimensional settings that $A^*ABB^*$ is EP. That is, $\mathcal R(A^*ABB^*) =\mathcal R((A^*ABB^*)^*).$ We could extend the result to infinite case, which upto our knowledge nobody has proved  for infinite dimensional case. 

\section{Preliminaries}

Let $\mathcal H_1$ and $\mathcal H_2 $ be Hilbert spaces and  $\mathcal B(\mathcal H_1, \mathcal H_2 )$ denote the set of all linear bounded operators from $\mathcal H_1$ to $\mathcal H_2$. We abbreviate $\mathcal B(\mathcal H_1)=\mathcal B(\mathcal H_1,\mathcal H_1)$. For $A \in \mathcal B(\mathcal H_1, \mathcal H_2 ),$ we denote by $A^*$, $\mathcal N(A)$ and $\mathcal R(A)$ respectively, the adjoint, the null-space and the range of $A$. An operator $A\in \mathcal B(\mathcal H_1)$ is said to be a projection if $A^2=A.$ A projection is said to be orthogonal if $A^2=A=A^*.$ The Moore-Penrose inverse of $A \in \mathcal B(\mathcal H_1, \mathcal H_2)$ is the
operator  $X \in \mathcal B(\mathcal H_2,\mathcal H_1 )$ which satisfies the Penrose equations
\begin{eqnarray}
AXA&=&A\label{eqn1}\\
XAX&=&X\label{eqn2}\\	
(AX)^* &=& AX\label{eq3}\\
(XA)^*&=& XA .\label{eqn4}
\end{eqnarray}

A matrix  $X$ is called a $\{i,\dots, j\}$-generalized inverse of $A,$ denoted by $A^{(i,\dots,j)}$ if it satisfies the $i^{th},\dots, j^{th}$ conditions of Penrose equations. The collection of all $\{i,\dots , j\}$-generalized inverses of $A$ is denoted by $A\{i,\dots,j\}$.
 If the Moore-Penrose inverse of $A$ exists, then it
is unique and it is denoted by $A^\dagger$. 

For the sake of clarity as well as for easier reference, we mention the following properties of Moore-Penrose inverse without proof \cite{Wangbook}. 
\begin{lemma}\label{properties}
	Let  $A \in \mathcal B(\mathcal H_1, \mathcal H_2)$ be a closed range
	operator. The following statements hold:
	\begin{enumerate}
		\item[(i)] $(A^\dagger)^\dagger=A.$
		\item[(ii)] $(A^\dagger)^*=(A^*)^\dagger.$ 
		\item[(iii)]   $A = AA^*(A^*)^\dagger = (A^*)^\dagger A^*A.$	
			\item[(iv)] $A^\dagger=A^*(AA^*)^\dagger=(A^*A)^\dagger A^*.$ 
		\item[(v)] $(AA^*)^\dagger =(A^*)^\dagger A^\dagger$, $(A^*A)^\dagger =A^\dagger(A^*)^\dagger .$
		\item[(vi)]	$A^* = A^*AA^\dagger = A^\dagger AA^*.$
		\item[(vii)] $\mathcal R(A) =\mathcal R(AA^*)= \mathcal R(AA^\dagger ).$
		\item[(viii)]  $\mathcal R(A^\dagger) =\mathcal R(A^*)= \mathcal R(A^\dagger A) =\mathcal R(A^*A).$
		\item[(ix)]  $AA^\dagger = P_{ \mathcal R(A )}$ and $ A^\dagger A =P_{\mathcal R(A^*)} =P_{\mathcal R(A^\dagger )}.$
		\item[(x)] $(A^n)\+=(A\+)^n$ for $n\ge 1.$
	\end{enumerate}
	Here,   $P_{\mathcal R(A )}$ and $P_{\mathcal R(A^* )}$ denote projection operators onto $\mathcal R(A )$ and $\mathcal R(A^*) $ respectively. We use $A{^\dagger}^* $ instead of $(A^\dagger)^*$  throughout the paper.
\end{lemma}

\begin{lemma}(\cite{Brock}, Theorem 1)\label{eplemma}
	Let $A\in \mathcal B(\mathcal H_1, \mathcal H_2)$ be a closed range operator such that  $\mathcal R(A)=\mathcal R(A^*).$ Then $AA^\dagger=A^\dagger A$ and $A^nA^\dagger=A^{n-1},n\ge 2.$
\end{lemma}
\begin{lemma} \label{projection_lemma} 
	Let $\mathcal H$ be a Hilbert space and $A\in \mathcal B(\mathcal H)$ be a projection. Then
	$A$ is Hermitian if and only if $A=AA^* A.$
\end{lemma}
\begin{proof}
	Suppose $A=AA^* A$ and $A$ is a projection. Let $ B=A-A^*$. Then it is easy to verify that  $B^3=0$ and $\mathcal R(B)=\mathcal R(B^*).$ By Lemma \ref{eplemma}, $B^3(B^\dagger )^2  =0$ gives $B=0.$
	Thus $A $ is Hermitian.	Converse follows directly. 
\end{proof}

\begin{lemma}\label{lemma134}(\cite{Xiong}, Lemma 1.3)
	Let $A\in \mathcal B(\mathcal H_1,\mathcal H_2)$ have a closed range and $B\in \mathcal B(\mathcal H_2,\mathcal H_1)$. Then
	\begin{enumerate}
		\item[(i)] $B \in A\{1, 3\} \Leftrightarrow A^*AB = A^*$
		\item[(ii)] $B \in A\{1, 4\} \Leftrightarrow BAA^* = A^*$.
	\end{enumerate}
\end{lemma} 
\begin{theorem}\label{Douglas}(\cite{Douglas}, Theorem 1)
	Let  $A$ and $B$ be bounded operators on a Hilbert space $\mathcal H.$ The following statements are equivalent:
\begin{enumerate}
	\item[(i)] $\mathcal R(A) \subseteq  \mathcal R(B);$
\item[(ii)] $ AA^* \le \lambda^2BB^*$ for some $\lambda\ge 0$ i.e., there exists $\lambda \ge 0$ such that $ \|A^*x\|\le \lambda \|B^*x\|$ for all $x \in \mathcal H;$
\item[(iii)]  there exists a bounded operator $C$ on $\mathcal H$ so that $A=BC.$
	\end{enumerate}
\end{theorem}
\begin{theorem}\label{selfadjoint}(\cite{hall}, Theorem 7.20)
Let  $A\in \mathcal B(\mathcal H)$ be self-adjoint. Then there exist a measure space $(X,\Sigma,\mu)$ and a bounded, measurable, real-valued function $f$ on $X$ and a unitary operator $U:\mathcal H\rightarrow L^2(X,\mu)$ such that $$A=U^*TU$$
where $T$ is the multiplication operator given by $(T\psi)(x)=f(x)\psi(x),x\in X.$ 
\end{theorem}
\begin{definition}
	Let $(\mathcal H,\langle.,.\rangle)$ be a Hilbert space and  $A\in \mathcal B(\mathcal H).$ $A$ is called a positive semi-definite operator if $\langle Ax,x\rangle\ge 0\text{  for all }  x \in \mathcal H.$
\end{definition}
\begin{lemma}\label{asquarea}
	Let $\mathcal H$ be a Hilbert space and  $A\in \mathcal B(\mathcal H)$ be a positive semi-definite operator such that $A^m=A^n$ for some natural numbers $m\ne n.$ Then $A^2=A.$
\end{lemma}
\begin{proof}
	We know that a positive semi-definite operator is self adjoint. By Theorem \ref{selfadjoint}, we can write $$A=U^*TU$$
	where $T$ is the multiplication operator given by $(T\psi)(x)=f(x)\psi(x),x\in X.$ 	Using the positive semi-definiteness of the operator we get $f(x)\ge0 \text{ } \forall x \in X.$
	
	It is given that $A^m=A^n$ which implies $f^m(x)\psi(x)=f^n(x)\psi(x) \text{ } \forall \psi(x)\in  L^2(X,\mu).$ In particular for $\psi(x)=1,$ we get $f^m(x)(1-f^{n-m}(x))=0,$ from which we can conclude $f(x)=0\text{ or } f(x)=1$ as $f(x)\ge 0.$

Now, $T^2\psi(x)=T(f(x)\psi(x))=f(x)^2\psi(x)=f(x)\psi(x)=T\psi(x)\text{ for all  } \psi(x)\in  L^2(X,\mu).$ Also, $U^*T^2U=U^*TU\Rightarrow A^2=A.$
 \end{proof}
\begin{lemma}\label{abba}
	Let $A$ and $B$ be orthogonal projections on a Hilbert space $\mathcal H$ and $m>n\ge 1.$ If $(ABA)^m=(ABA)^n,$ then $AB=BA.$
\end{lemma}
\begin{proof}
	$ABA=ABBA=ABB^*A^*=AB(AB)^*.$ Thus $ABA$ is Hermitian and positive semi-definite as $AB(AB)^*$ is so. Then by Lemma \ref{asquarea},  $(ABA)^m=(ABA)^n$ implies  $(ABA)^2=ABA.$ Consider $(ABA-AB)(ABA-AB)^*=(ABA-AB)(ABA-BA)=(ABA)^2-ABABA-ABABA+ABA=0.$ Thus $ABA=AB.$ Similarly, we can verify $(ABA-BA)(ABA-BA)^*=0$ which gives $ABA=BA.$ Thus, we get $AB=BA.$
\end{proof}

\begin{lemma}\label{aBBA}
	Let $A$ and $B$ be orthogonal projections on a Hilbert space $\mathcal H$ and $m>n\ge1.$ If $(AB)^m=(AB)^n,$ then $AB=BA.$
\end{lemma} 
\begin{proof}
	Since $(AB)^2A=ABABA=ABA ABA=(ABA)^2,$ thus $(AB)^mA=(ABA)^m$ for all $m\ge1.$ Now it is clear that$(AB)^m=(AB)^n$ gives $(ABA)^m=(ABA)^n.$
 Then by Lemma \ref{abba}, we get $AB=BA.$ 
\end{proof}

\section{Main Results}
We start the section with some examples to show that reverse order law does not  hold good for closed range Hilbert space operators in general.  
\begin{example}\label{ex1}
	Let $\mathcal H=\ell_2$ be the space of all square summable sequences.  For $x=(x_1,x_2,x_3,\dots)\in\mathcal H,$ define $A(x)=(x_1+x_2,x_2,x_3,x_4,\dots)$ and $B(x)= (x_1,0,x_3,0,x_5,\dots).$ Then 	$AB(x)=A(x_1,0,x_3,0,x_5,\dots)=(x_1,0,x_3,0,x_5,\dots)=Bx$. It can be  verified easily that $A,B$ and $AB$ are bounded and have closed ranges. We see  that $A^*(x)=(x_1,x_1+x_2,x_3,x_4,\dots)$ and $B^*(x) =(x_1,0,x_3,0,x_5,\dots)=Bx.$ Using computational methods for Moore-Penrose inverse of operators (\cite{Wangbook}, p.327) we get
	\begin{equation*}
	A^\dagger(x)=(x_1-x_2,x_2,x_3,x_4,\dots),\text{ } B^\dagger =B \quad and \quad (AB)^\dagger=B^\dagger.
	\end{equation*}
	Hence $B^\dagger A^\dagger(x)=B^\dagger(x_1-x_2,x_2,x_3,x_4,\dots)=(x_1-x_2,0,x_3,0,x_5,\dots) \ne (AB)^\dagger(x),$ thus $(AB)\+\ne B\+A\+.$
\end{example}	

\begin{example}\label{ex2}
	Let $\mathcal H=\ell_2.$  For $x=(x_1,x_2,x_3,\dots)\in \mathcal H,$ define $A(x)=(0,x_2,0,x_4,0,\dots) $ and $B(x)= (x_1+x_2,2x_1+2x_2,x_3,x_4,\dots).$ Then $ABx=(0,,2x_1+2x_2,0,x_4,\dots).$ It is easy to verify that  $A,B$ and $AB$ are bounded and have closed ranges.
	Also  $A^*(x)= (0,x_2,0,x_4,0,x_6,\dots)=A^\dagger x,$ and $B^*(x) =(x_1+2x_2,x_1+2x_2,x_3,x_4,\dots), B^\dagger x=(\frac{1}{10}(x_1+2x_2),\frac{1}{10}(x_1+2x_2),x_3,x_4,\dots).$ Thus we get	$$B^\dagger A^\dagger(x)=B^\dagger (0,x_2,0,x_4,0,x_6,\dots)=(\frac{x_2}{5},\frac{x_2}{5},0,x_4,0,\dots)$$ and $$(AB)^\dagger x=\\(\frac{x_2}{4},\frac{x_2}{4},0,x_4,0,\dots).	$$ Hence $(AB)^\dagger \ne B^\dagger A^\dagger .$ One can also check that $B^\dagger A^\dagger$ satisfies the third and fourth but not the first and second Penrose equations.	
\end{example}		

\begin{lemma}(\cite{Izumino}, Proposition 2.1) \label{closedrangelemma}
	Let $\mathcal H_1,\mathcal H_2,\mathcal H_3$ be Hilbert spaces, and let $A\in \mathcal B(\mathcal H_2,\mathcal H_3)$ and $B\in \mathcal B(\mathcal H_1,\mathcal H_2)$ be such that $A, B$ have closed ranges. Then $AB$ has a closed range if and only if $A\+ABB\+$ has a closed range.
\end{lemma}
The  results mentioned below  in Theorems  \ref{equivalence} to \ref{equivalence1234} are proved in  $C^*$-algebra setting \cite{Dijana}. For the sake of completeness, we give the proof of those in Hilbert space setting. However, our proofs are simpler than  those  available for reverse order law for closed range Hilbert space operators.

In the following result,  the existence of  $(A^\dagger ABB^\dagger )^\dagger$ is guaranteed by Lemma \ref{closedrangelemma}.
\begin{theorem} \label{equivalence}
	Let $\mathcal H_1,\mathcal H_2,\mathcal H_3$ be Hilbert spaces, and let $A\in \mathcal B(\mathcal H_2,\mathcal H_3)$ and $B\in \mathcal B(\mathcal H_1,\mathcal H_2)$ be such that $A, B, AB$ have closed ranges.  Then the following statements are equivalent:
	\begin{enumerate}
		\item[(i)] $ABB^\dagger A^\dagger AB=AB$ ;
		\item[(ii)] $B^\dagger A^\dagger ABB^\dagger A^\dagger =B^\dagger A^\dagger $ ;
		\item[(iii)] $BB^\dagger A^\dagger A$ is a projection ;
	\item[(iv)]	$A^\dagger ABB^\dagger =BB^\dagger A^\dagger A$ ;
				\item[(v)]  $A^\dagger ABB^\dagger $ is a projection ;
		\item[(vi)] $(A^\dagger ABB^\dagger )^\dagger =BB^\dagger A^\dagger A$;
		\item[(vii)] $B^\dagger (A^\dagger ABB^\dagger )^\dagger A^\dagger =B^\dagger A^\dagger .$ 
	\end{enumerate}	
\end{theorem}

\begin{proof}
	\noindent \underline{$(i) \Rightarrow (ii)$}:
	If $ABB^\dagger A^\dagger AB=AB$, then
	\begin{align*}
	B^\dagger A^\dagger &=(B^*B)^\dagger B^*A^*(AA^*)^\dagger \\
	&=(B^*B)^\dagger (AB)^*(AA^*)^\dagger \\
	&=(B^*B)^\dagger (ABB^\dagger A^\dagger AB)^*(AA^*)^\dagger \\
	&=(B^*B)^\dagger B^*A^\dagger ABB^\dagger A^*(AA^*)^\dagger \\
	&=B^\dagger A^\dagger ABB^\dagger A^\dagger . 
	\end{align*}
	\underline{ $(ii)\Rightarrow(iii)$}:
	Using $(ii)$ we see that   $(BB^\dagger A^\dagger A)^2=BB^\dagger A^\dagger ABB^\dagger A^\dagger A=BB^\dagger A^\dagger A$. Hence it shows that $BB^\dagger A^\dagger A$ is a projection. \\
		\underline{ $(iii)\Rightarrow(iv)$}: Consider	\begin{align*}BB^\dagger A^\dagger A (BB^\dagger A^\dagger A)^* BB^\dagger A^\dagger A
	&=BB^\dagger A^\dagger A (A^\dagger A)^* (BB^\dagger)^* BB^\dagger A^\dagger A\\
	&=BB^\dagger A^\dagger A(A^\dagger A) (BB^\dagger)BB^\dagger A^\dagger A\\
	&=BB^\dagger A^\dagger A BB^\dagger A^\dagger A\\
	&=BB^\dagger A^\dagger A.
	\end{align*}
	Then by Lemma \ref{projection_lemma}, we get 	$(BB^\dagger A^\dagger A)^*=BB^\dagger A^\dagger A.$ Thus
	$BB^\dagger A^\dagger A=A^\dagger ABB^\dagger $.\\
	\underline{ $(iv)\Rightarrow (v)$}: It is given that $A^\dagger ABB^\dagger =BB^\dagger A^\dagger A$. We have
	\begin{align*}
	(A^\dagger ABB^\dagger )^2&=A^\dagger ABB^\dagger A^\dagger ABB^\dagger \\
	&=BB^\dagger A^\dagger A A^\dagger ABB^\dagger\\
	&=BB^\dagger A^\dagger A BB^\dagger\\
	&=A^\dagger ABB^\dagger BB^\dagger \\
	&=A^\dagger A BB^\dagger.
	\end{align*}
	\underline{$(v) \Rightarrow (vi)$}: Using the fact that $A$ is a projection if and only if $A^*$ is a projection, it is easy to verify   all Penrose equations.\\
	\underline{$(vi)\Rightarrow(vii)$}: Pre- and post-multiplying by $B^\dagger$ and $A^\dagger$ respectively in $(vi)$, we get the desired result.\\
	\underline{$(vii)\Rightarrow (i)$}: We have
	\begin{align*}
ABB^\dagger A^\dagger AB&=ABB^\dagger (A^\dagger ABB^\dagger )^\dagger A^\dagger AB\\&=AA^\dagger ABB^\dagger (A^\dagger ABB^\dagger )^\dagger A^\dagger ABB^\dagger B\\&=AA^\dagger ABB^\dagger B\\&=AB.
	\end{align*}\end{proof}
Next, we give nine equivalent conditions for $B\+A\+$ to be a $\{1,2,3\}$-generalized inverse of $AB$ in Hilbert space settings. The existence of  $(ABB^\dagger )\+$ follows as the  ranges of  $ABB^\dagger $ and $AB$ are equal.
\begin{theorem}\label{equivalence123}
	Let $\mathcal H_1,\mathcal H_2,\mathcal H_3$ be Hilbert spaces, and let $A\in\mathcal  B(\mathcal H_2,\mathcal H_3)$ and $B\in \mathcal B(\mathcal H_1,\mathcal H_2)$ be such that $A, B, AB$ have closed ranges.  Then the following statements are equivalent:
	\begin{enumerate}
		\item[(i)]$AB(AB)^\dagger =ABB^\dagger A^\dagger ;$
		\item[(ii)] $B^\dagger A^\dagger \in AB\{1,2,3\};$
		\item[(iii)] $BB^\dagger A^*AB=A^*AB;$
		\item[(iv)] $(AB)(AB)^\dagger A=ABB^\dagger ;$
		\item[(v)] $A^*ABB^\dagger =BB^\dagger A^*A;$
		\item[(vi)] $(ABB^\dagger )^\dagger =BB^\dagger A^\dagger ;$
		\item[(vii)] $B^\dagger (ABB^\dagger )^\dagger =B^\dagger A^\dagger ;$
		\item[(viii)] $B\{1,3\}A\{1,3\} \in AB\{1,3\};$
		\item[(ix)] $B^\dagger A^\dagger \in AB\{1,3\}.$
	\end{enumerate}	
\end{theorem}
\begin{proof} We prove the equivalency of  all the statements in the following order of  implications:\\ $(i)\Rightarrow(ii)\Rightarrow(iii)\Rightarrow(iv)\Rightarrow(v)\Rightarrow(vi)\Rightarrow(vii)\Rightarrow(ix)\Leftrightarrow(viii),(ix)\Rightarrow(i).$\\
	\underline{$(i)\Rightarrow(ii)$}: Since $ABB^\dagger A^\dagger =AB(AB)^\dagger ,$ post-multiplying by $AB$ we get,
	$ABB^\dagger A^\dagger AB=AB(AB)^\dagger AB=AB$.
	Hence $B^\dagger A^\dagger \in AB\{1\}$. By Theorem \ref{equivalence}, $B^\dagger A^\dagger \in AB\{2\}$.
	Now, $(ABB^\dagger A^\dagger )^*=(AB(AB)^\dagger )^*
	=AB(AB)^\dagger 
	=ABB^\dagger A^\dagger.
	$
	Thus $B^\dagger A^\dagger \in AB\{1,2,3\}$.\\  	
	\underline{$(ii) \Rightarrow(iii)$}:
	Suppose $B^\dagger A^\dagger \in AB\{1,2,3\}$. Then by Theorem \ref{equivalence}, 
	$A^\dagger ABB^\dagger =BB^\dagger A^\dagger A$.
	Thus we get
	\begin{align*}
	A^*AB&=A^*AA^\dagger ABB^\dagger B\\
	&=A^*ABB^\dagger A^\dagger AB\\
	&=A^*(ABB^\dagger A^\dagger )^*AB\\
	&=A^*{A^\dagger }^*BB^\dagger A^*AB\\
	&=(A^\dagger A)^*BB^\dagger A^*AB\\
	&=A^\dagger ABB^\dagger A^*AB\\
	&=BB^\dagger A^\dagger AA^*AB\\
	&=BB^\dagger A^*AB.
	\end{align*}
	\underline{$(iii)  \Rightarrow (iv)$}: $(AB)(AB)^\dagger A=(A^*AB(AB)^\dagger)^*=(BB^\dagger A^*AB(AB)^\dagger)^*=AB(AB)^\dagger ABB^\dagger= ABB^\dagger.$\\
	\underline{$(iv)\Rightarrow(v)$}: $A^*ABB^\dagger=A^*AB(AB)^\dagger A.$ Hence $A^*ABB^\dagger$ is Hermitian and $A^*ABB^\dagger=BB^\dagger A^*A.$\\
	\underline{$(v)\Rightarrow(vi)$}: We show this by verifying all Penrose equations. Given that $A^*ABB^\dagger=BB^\dagger A^*A.$
 Pre-multiplying by ${A^\dagger}^*,$ we get $ABB^\dagger={A^\dagger}^*BB^\dagger A^*A.$ Hence $(ABB^\dagger)(BB^\dagger A^\dagger)(ABB^\dagger)=ABB^\dagger A^\dagger ABB^\dagger= {A^\dagger}^* BB^\dagger A^*AA^\dagger ABB^\dagger={A^\dagger}^*A^*ABB^\dagger BB^\dagger= ABB^\dagger.$ This shows that $BB^\dagger A^\dagger \in ABB^\dagger \{1\}. 
 $ Now $BB^\dagger A^\dagger=BB^\dagger A^* (AA^*)^\dagger =(ABB^\dagger)^*(AA^*)^\dagger=(ABB^\dagger B B^\dagger A^\dagger ABB^\dagger)^*(AA^*)^\dagger=(ABB^\dagger A^\dagger ABB^\dagger)^*(AA^*)^\dagger=BB^\dagger A^\dagger ABB^\dagger A^*(AA^*)^\dagger=BB^\dagger A^\dagger ABB^\dagger BB^\dagger A^\dagger.$ 
 Thus $BB^\dagger A^\dagger \in ABB^\dagger\{1,2\}.$
 Also $(ABB^\dagger)(BB^\dagger A^\dagger)=({A^\dagger}^*A^*ABB\+) (BB^\dagger A^\dagger)={A^\dagger}^*BB^\dagger A^*ABB^\dagger A^\dagger$
and $(BB^\dagger A^\dagger)(ABB^\dagger)$  are Hermitian. It ensures that $(ABB^\dagger)^\dagger=BB^\dagger A^\dagger.$\\
	\underline{$(vi)\Rightarrow(vii)$}: Pre-multiplying the given condition  by $B^\dagger $, we get $(ii)$. \\
\underline{$(vii)\Rightarrow(ix)$}: It is clear that $ABB^\dagger (ABB^\dagger)^\dagger ABB^\dagger =ABB^\dagger.$ Then by $(vii),$ we have $ABB^\dagger A^\dagger ABB^\dagger =ABB^\dagger .$ Post-multiplying by $B$ we get $ABB^\dagger A^\dagger AB.$ Thus $B^\dagger A^\dagger \in AB\{1\}.$ Also, $ABB^\dagger (ABB^\dagger)^\dagger =ABB^\dagger A^\dagger$ is Hermitian. Thus $B^\dagger A^\dagger \in AB\{1,3\}.$\\	\underline{$(ix)\Rightarrow(viii)$}: Let $CD\in B\{1,3\}A\{1,3\}$ where $C\in B\{1, 3\}$ and $D\in A\{1, 3\}$. By Lemma \ref{lemma134}, $ C$ and $D$ satisfy $B^*BC=B^*$ and $A^*AD = A^*$ as $B^\dagger BC=B^\dagger$ and $ A^\dagger A D=A^\dagger$. Since $B^\dagger A^\dagger \in AB\{1,3\}$ and using Theorem \ref{equivalence} $(iv),$ we get
	\begin{align*}
	(AB)^*(AB)CD&=(ABB^\dagger A^\dagger AB)^*ABCD\\
	&=(AB)^*A A^\dagger ABB^\dagger BCD\\
	&=(AB)^*A A^\dagger ABB^\dagger D\\
	&=(AB)^*A BB^\dagger A^\dagger A D\\
	&=(AB)^*A BB^\dagger A^\dagger\\ 
	&=(ABB^\dagger A^\dagger AB)^*\\ 
	&= (AB)^*.
	\end{align*}
	\underline{$(viii)\Rightarrow(ix)$}: Obvious.	\\
	\underline{$(ix)\Rightarrow(i)$}: By assumption, we have $ABB^\dagger A^\dagger AB =AB $ and $ (ABB^\dagger A^\dagger)^*=ABB^\dagger A^\dagger.$ Post-multiplying by $(AB)^\dagger$ in the first equation and taking adjoint on both sides, we get 
$AB(AB)^\dagger =ABB^\dagger A^\dagger .$
\end{proof}
The following result is similar to Theorem \ref{equivalence123}. It gives nine equivalent conditions for $B^\dagger A^\dagger$ to be a $\{1,2,4\}$-generalized inverse of $AB$ in Hilbert space setting. Here the existence of $ (A^\dagger AB)^\dagger$ is guaranteed as the ranges of $ (A^\dagger AB)^*$ and $ (AB)^* $ are the same. 
\begin{theorem}\label{equivalence124}	
	Let $ \mathcal H_1, \mathcal  H_2, \mathcal  H_3$ be Hilbert spaces, and let $A\in\mathcal  B( \mathcal H_2, \mathcal  H_3)$ and $B\in \mathcal B( \mathcal H_1, \mathcal  H_2)$ be such that $A, B, AB$ have closed ranges.  Then the following statements are equivalent:
	\begin{enumerate}
		\item[(i)]$(AB)^\dagger AB=B^\dagger A^\dagger AB;$
		\item[(ii)] $B^\dagger A^\dagger \in AB\{1,2,4\};$
		\item[(iii)] $ABB^*=ABB^*A^\dagger A;$
		\item[(vi)] $B(AB)^\dagger AB=A^\dagger AB;$
		\item[(v)] $A^\dagger ABB^*=BB^*A^\dagger A;$
		\item[(vi)] $(A^\dagger AB)^\dagger =B^\dagger A^\dagger A;$
		\item[(vii)] $(A^\dagger AB)^\dagger A^\dagger =B^\dagger A^\dagger;$
		\item[(viii)] $B\{1,4\}A\{1,4\} \in AB\{1,4\};$
		\item[(ix)] $B^\dagger A^\dagger \in AB\{1,4\}.$
	\end{enumerate}
\end{theorem}
\begin{proof}
Proof is similar to Theorem \ref{equivalence123}.
\end{proof}
\begin{theorem}  \label{equivalence1234}	 
	Let $ \mathcal H_1, \mathcal  H_2, \mathcal  H_3$ be Hilbert spaces, and let $A\in\mathcal  B( \mathcal H_2, \mathcal  H_3)$ and $B\in \mathcal B( \mathcal H_1,  \mathcal H_2)$ be such that $A, B, AB$ have closed ranges.  Then the following statements are equivalent:
	\begin{enumerate}
		\item[(i)]$(AB)^\dagger =B^\dagger A^\dagger ;$
		\item[(ii)] $(AB)(AB)^\dagger =ABB^\dagger A^\dagger $ and $(AB)^\dagger AB=B^\dagger A^\dagger AB;$
		\item[(iii)] $A^*AB=BB^\dagger A^*AB$ and $ABB^*=ABB^*A^\dagger A;$
		\item[(iv)] $AB(AB)^\dagger A=ABB^\dagger $ and $B(AB)^\dagger AB=A^\dagger AB;$
		\item[(v)] $A^*ABB^\dagger =BB^\dagger A^*A$  and $BB^*A^\dagger A=A^\dagger ABB^*;$ 
		\item[(vi)] $(ABB^\dagger)^\dagger=BB^\dagger A^\dagger \text{and}~ (A^\dagger AB)^\dagger =B^\dagger A^\dagger A;$
		\item[(vii)] $B^\dagger (ABB^\dagger)^\dagger=B^\dagger A^\dagger \text{and}~ (A^\dagger AB)^\dagger A^\dagger=B^\dagger A^\dagger;$
		\item[(viii)] $B\{1,3\}A\{1,3\} \in AB\{1,3\}$ and $B\{1,4\}A\{1,4\} \in AB\{1,4\};$
		\item[(ix)] $B^\dagger A^\dagger \in AB\{1,3,4\}.$
	\end{enumerate}	 
\end{theorem} 
\begin{proof}
Follows from Thorems \ref{equivalence123} and \ref{equivalence124}.
\end{proof}
\begin{remark}
	Consider the operators $A$ and $B$ on $ \mathcal H$ defined in Example \ref{ex1}. Then for all $x\in  \mathcal H,$ $A^*ABx=(x_1,x_1,x_3,0,x_5,\dots)$ and $
	BB^\dagger A^*ABx=(x_1,0,x_3,0,x_5,\dots).$ Hence $A^*AB\ne BB^\dagger A^*AB$ and $ABB^*x=(x_1,0,x_3,0,x_5,\dots)=ABB^*A^\dagger Ax.$ Note that the conditions in $(iii)$ of Theorem \ref{equivalence1234} are not satisfied and $(AB)\+\ne B\+A\+$   which was shown in Example \ref{ex1}.
\end{remark}

\begin{theorem} \label{rangethm}  
	Let $  \mathcal H_1, \mathcal  H_2, \mathcal  H_3$ be Hilbert spaces, and let $A\in\mathcal  B( \mathcal H_2, \mathcal  H_3)$ and $B\in \mathcal B( \mathcal H_1, \mathcal  H_2)$ be such that $A, B, AB$ have closed ranges. Then the following statements hold :
	\begin{enumerate}
		\item [(i)] $B^\dagger =(AB)^\dagger A \Leftrightarrow \mathcal R(B)=\mathcal R(A^* AB).$
		\item [(ii)] $A^\dagger=B(AB)^\dagger \Leftrightarrow \mathcal R(A^*)= \mathcal R(BB^* A^*).$
	\end{enumerate}
\end{theorem}
\begin{proof}
	$(i)$ Suppose that $B^\dagger =(AB)^\dagger A.$ Pre-multiplying by $AB,$ we get
	$ABB^\dagger=(AB) (AB)^\dagger A,$ which is equivalent to 			$BB^\dagger A^*AB=A^*AB$ by Theorem \ref{equivalence123}. This implies
	$\mathcal R(A^* AB) \subseteq \mathcal R(B).$ Now $B^\dagger=(AB)^\dagger A=[(AB)^* (AB)]^\dagger (AB)^* A $ implies that $B^{\dagger^*}=A^* AB [(AB)^* (AB)]^\dagger.$
	Thus \begin{equation*}
	\mathcal R(B)=\mathcal R(B^{\dagger ^*})=\mathcal R(A^* AB[(AB)^* (AB)]^\dagger)\subseteq \mathcal R(A^* AB).
	\end{equation*} 
	Conversely, $\mathcal R(B)=\mathcal R(A^* AB)\implies BB^\dagger A^*AB=A^*AB.$ By Theorem \ref{equivalence123}, $B^\dagger A^\dagger \in AB\{1,2,3\},$ $\mathcal R({B^\dagger}^*)=\mathcal R(B)=\mathcal R(A^* AB) \subseteq \mathcal R(A^*)$ gives   $A^\dagger AB=B$ and $A^\dagger A B^{\dagger^ \star}=B^{\dagger^ \star}\implies B^\dagger A^\dagger A=B^\dagger.$ It shows that $(A^\dagger AB)^\dagger =B^\dagger =B^\dagger A^\dagger A.$
By Theorem \ref{equivalence124}, $B^\dagger A^\dagger \in AB\{1,2,4\}$ and hence $ (AB)^\dagger =B^\dagger A^\dagger.$ Therefore $B^\dagger=B^\dagger A^\dagger A=(AB)^\dagger A.$\\
	$(ii)$ Proof is similar to $(i).$
\end{proof}

\begin{theorem}\label{rangeinclusion} 	
	Let $ \mathcal H_1, \mathcal  H_2, \mathcal  H_3$ be Hilbert spaces, and let $A\in\mathcal  B( \mathcal H_2, \mathcal  H_3)$ and $B\in \mathcal B( \mathcal H_1, \mathcal  H_2)$ be such that $A, B, AB$ have closed ranges.   Then the following statements hold :
	
	\begin{enumerate}
		\item [(i)] $(AB)^\dagger=(A^\dagger AB)^\dagger A^\dagger \Leftrightarrow \mathcal R(AA^* AB)=\mathcal R(AB).$
		\item [(ii)] $(AB)^\dagger=B^\dagger (ABB^\dagger)^\dagger\Leftrightarrow \mathcal R(B^* B(AB)^*)=\mathcal R((AB)^*).$
	\end{enumerate}
\end{theorem}
\begin{proof}	$(i)$ If we replace $A$ by $A\+$ and $B$ by $AB$ in Theorem \ref{rangethm} $(i),$ we get $(AB)^\dagger=(A^\dagger AB)^\dagger A^\dagger \Leftrightarrow \mathcal R(AB)=\mathcal R(A{\+}^*A\+AB)=\mathcal R(A{\+}^*B).$ Now by Theorem \ref{Douglas}, there exists a bounded operator $C$ such that $AB= A{\+}^*BC.$ Pre-multiplying by $AA^*$ we get $AA^*AB= AA^*A{\+}^*BC=A(A\+A)^*BC=ABC.$ Thus we get $\mathcal R(AA^* AB)\subseteq \mathcal R(AB).$ Similarly, we can prove $\mathcal R(AB)\subseteq \mathcal R(AA^* AB)   .$
	
		$(ii)$ Replace $A$ by $AB$ and $B$ by $B\+$ in Theorem \ref{rangethm} $(ii)$ and use a similar argument as above.
\end{proof}

\begin{theorem}\label{Rolrange} Let $ \mathcal H_1,  \mathcal H_2, \mathcal  H_3$ be Hilbert spaces, and let $A\in\mathcal  B( \mathcal H_2,  \mathcal H_3)$ and $B\in \mathcal B( \mathcal H_1, \mathcal  H_2)$ be such that $A, B, AB$ have closed ranges.  Then the following statements are equivalent:
	\begin{enumerate}
		\item[(i)]	$(AB)^\dagger=B^\dagger A^\dagger ;$
		\item[(ii)] $\mathcal R(A^*AB)\subseteq \mathcal R(B)$ and $\mathcal R(BB^*A^*)\subseteq \mathcal R(A^*);$
		\item[(iii)]$\mathcal R(AA^*AB) = \mathcal R( AB)$ and $\mathcal R(BB^*A^*) \subseteq \mathcal R(A^*);$
		\item[(iv)] $\mathcal R(A^*AB)\subseteq \mathcal R(B)$ and $\mathcal R[(ABB^*B)^*] = \mathcal R[( AB)^*];$
	\item[(v)] $ \mathcal R(A^*ABB^*)=\mathcal R(BB^*A^*A).$ 
	\end{enumerate}
\end{theorem}
\begin{proof}\underline{$(i)\Leftrightarrow(ii)$}: Follows from Theorem \ref{equivalence1234} $(iii)$.\\
\underline{$(ii)\Rightarrow(iii)$}: By Theorem \ref{equivalence1234} $(vii),$ $(AB)\+=B\+A\+=(A\+AB)\+A\+.$ Then $(iii)$ follows from Threorem \ref{rangeinclusion} $(i).$\\
\underline{$(iii)\Rightarrow(v)$}: By Theorem \ref{Douglas}, there exists an operator $T$ such that $AA^*AB=ABT.$ Pre-multiplying by $A\+$ we get $A^*AB=A\+ABT=A\+ABB^*B{\+}^*T=BB^*A\+AB{\+}^*T
$ by Theorem \ref{equivalence124} $(v).$ Also, $A^*ABB^*=BB^*A\+AB{\+}^*TB^*. $ Thus 	 $ \mathcal R(A^*ABB^*)\subseteq \mathcal R(BB^*A^*A).$
Similarly, $AB=AA^*ABS,$ for some operator $S.$ Pre-multiplying by $A\+$ and post-multiplying by $B^*A^*A,$ we get $A\+ABB^*A^*A=A^*ABSB^*A^*A.$ Then by Theorem \ref{equivalence124} $(v),$  $BB^*A\+AA^*A=A^*ABSB^*A^*A.$ Thus $BB^*A^*A=A^*ABSB^*A^*A.$ It shows that $ \mathcal R(A^*ABB^*)= \mathcal R(BB^*A^*A).$\\
\underline{$(v)\Rightarrow(ii)$}:	By Theorem \ref{Douglas}, there exists an operator $T$ such that $A^*ABB^*=BB^*A^*AT.$ 
Pre-multiplying by $BB\+$ and post-multiplying by $B{\+}^*$ we get, $ BB\+A^*ABB^*B{\+}^*=BB^*A^*ATB{\+}^*.$ Hence $BB\+A^*AB=BB^*A^*ATB{\+}^*=A^*ABB^*B{\+}^*=A^*AB,$  thus $ \mathcal R(A^*AB)\subseteq \mathcal R(B).$ Similarly, we can prove $A\+ABB^*A^*=BB^*A^*$ and hence $\mathcal R(BB^*A^*)\subseteq \mathcal R(A^*).$ 
\end{proof}
\begin{remark}
	Let $A$ and $B$ be as defined in Example \ref{ex1}. Then $A^*ABx=(x_1,x_1,x_3,0,x_5,...).$ Thus $ \mathcal R(A^*AB)\not\subseteq \mathcal R(B),$ 
	$A^*(x)=(x_1,x_1+x_2,x_3,x_4,...)$ and $BB^*A^*x=(x_1,0,x_3,0,x_5,...)$ implies $\mathcal R(BB^*A^*)\subseteq \mathcal R(A^*).$ This shows that $ \mathcal R(A^*AB)\subseteq \mathcal R(B)$ is indispensable for reverse order law to hold. 
\end{remark}

By Lemma \ref{closedrangelemma}, $\mathcal R(AB) $ is closed if and only if $\mathcal R(A\+ABB\+)=\mathcal R(A^*A\+{^*}BB\+)=\mathcal R(A\+AB\+{^*}B^*)$ is closed. This happens if and only if $\mathcal R(A\+{^*}B)$ and $\mathcal R(AB\+{^*})$ are closed. Thus $(A\+{^*}B)\+$ and $(AB\+{^*})\+$ exist. Also, $\mathcal R(A)$ is closed if and only if $\mathcal R(A^*)$ is closed implies $\mathcal R(BB\+A\+A)$ is closed and hence $\mathcal R(B\+A\+)$ is closed. For natural numbers $m$ and $n,$ the existence of Moore-Penrose inverse of $(AA^*)^m$ and $(B^*B)^n$ is guaranteed as they are powers of Hermitian operators with closed ranges, according to the spectral mapping theorem. The existence of the Moore-Penrose inverse of all other operators discussed below can be guaranteed with the closedness of the ranges of $AB,A\+{^*}B,AB\+{^*}$ and $ B\+A\+. $

\begin{theorem} Let $ \mathcal H_1, \mathcal  H_2,  \mathcal H_3$ be Hilbert spaces, and let $A\in\mathcal  B( \mathcal H_2, \mathcal  H_3)$ and $B\in \mathcal B( \mathcal H_1, \mathcal  H_2)$ be such that $A, B, AB$ have closed ranges.  Then the following statements are equivalent:
	\begin{enumerate}
		\item[(1)]	$(AB)^\dagger=B^\dagger A^\dagger ;$
		\item[(2)] $B(AB)^\dagger A=BB^\dagger A^\dagger A;$
		\item[(3)] $AA^*(B^*A^*)^\dagger B^*B=AB;$
		\item[(4)] $(AB)^\dagger=B^\dagger A^\dagger ABB^\dagger A^\dagger ;$
		\item[(5)] 	$(AB)^\dagger=(A^\dagger AB)^\dagger A^\dagger$ and $(A^\dagger AB)^\dagger =B^\dagger A^\dagger A;$
		\item[(6)] $(AB)^\dagger=B^\dagger(ABB^\dagger)^\dagger$ and $(ABB^\dagger)^\dagger =BB^\dagger A^\dagger;$
		\item[(7)] $(AB)^\dagger=B^\dagger (A^\dagger ABB^\dagger )^\dagger A^\dagger$ and $(A^\dagger ABB^\dagger )^\dagger=BB^\dagger A^\dagger A ;$
		\item[(8)] $B^\dagger A^\dagger \in AB\{1,3,4\};$
		\item[(9)] $(AB)(AB)^\dagger =ABB^\dagger A^\dagger =A^\dagger{^*}BB^\dagger A^*$ and $(AB)^\dagger AB=B^\dagger A^\dagger AB=B^*A^\dagger A(B^\dagger)^*;$
		\item[(10)] $(A^\dagger{^*}B)^\dagger=B^\dagger A^* ;$
		\item[(11)]$A^\dagger (B^* A^\dagger)^\dagger B^*=A^\dagger ABB^\dagger;$
		\item[(12)]$AA^\dagger(B^*A^\dagger)^\dagger B^*B=AB;$
		\item[(13)]$(A^\dagger{^*}B)^\dagger=B^\dagger A^\dagger ABB^\dagger A^* ;$
		\item[(14)] $(A\+{^*}B)^\dagger=(A^\dagger AB)^\dagger A^*$ and $(A^\dagger AB)^\dagger =B^\dagger A^\dagger A;$
		\item[(15)]$(A\+{^*}B)^\dagger=B^\dagger(A\+{^*}BB^\dagger)^\dagger$ and $(A\+{^*}BB^\dagger)^\dagger =BB^\dagger A^*;$
		\item[(16)] $(A\+{^*}B)^\dagger=B^\dagger (A^\dagger ABB^\dagger )^\dagger A^*$ and $(A^\dagger ABB^\dagger )^\dagger=BB^\dagger A^\dagger A ;$
		\item[(17)]$B^\dagger A^*\in A\+{^*}B\{1,3,4\};$
		\item[(18)] $(B^*A\+)\+B^*A\+=ABB\+A\+=A\+{^*}BB\+A^*$ and $ B^*A\+(B^*A\+)\+=B\+A\+AB=B^*A\+AB\+{^*};$
		\item[(19)]$(AB\+{^*})^\dagger=B^* A^\dagger ;$
		\item[(20)]$B\+{^*}(AB\+{^*})^\dagger A=BB^\dagger A^\dagger A;$
		\item[(21)]$AA^*(B\+A^*)^\dagger B\+B=AB;$
		\item[(22)] $(AB\+{^*})^\dagger=B^* A^\dagger ABB^\dagger A^\dagger ;$
		\item[(23)]	$(AB\+{^*})^\dagger=(A^\dagger AB\+{^*})^\dagger A^\dagger$ and $(A^\dagger AB\+{^*})^\dagger =B^*A^\dagger A;$
		\item[(24)]$(AB\+{^*})^\dagger=B^*(ABB^\dagger)^\dagger$ and $(ABB^\dagger)^\dagger =BB^\dagger A^\dagger;$
		\item[(25)] $(AB\+{^*})^\dagger=B^* (A^\dagger ABB^\dagger )^\dagger A^\dagger$ and $(A^\dagger ABB^\dagger )^\dagger=BB^\dagger A^\dagger A ;$
		\item[(26)] $B^* A^\dagger \in AB\+{^*}\{1,3,4\};$
		\item[(27)]$(B\+A^*)\+B\+A^*=ABB\+A\+=A\+{^*}BB\+A^*$ and $ B\+A^*(B\+A^*)\+=B\+A\+AB=B^*A\+AB\+{^*};$
		\item[(28)] $(B\+A\+)\+=AB;$
		\item[(29)]$A\+(B\+A\+)\+B\+=A\+ABB\+;$
		\item[(30)]$(AA^*)\+(B\+A\+)\+(B^*B)\+=A\+{^*}A\+{^*};$
		\item[(31)]$(B\+A\+)\+=ABB\+A\+AB;$	
	    \item[(32)]$(B\+A\+)\+=A(B\+A\+A)\+$ and $(B\+A\+A)\+=A\+AB;$			
		\item[(33)]$(B\+A\+)\+=(BB\+A\+)\+B$ and $(BB\+A\+)\+=ABB\+;$		
		\item[(34)]		$(B\+A\+)\+=A(BB\+A\+A)\+B$ and $ (BB\+A\+A)\+=A\+ABB\+;$
		\item[(35)]		$AB \in B\+A\+\{1,3,4\};$
		\item[(36)]	$	B\+A\+(B\+A\+)\+=B\+A\+AB=B^*A\+AB\+{^*}$ and $(B\+A\+)\+B\+A\+=ABB\+A\+=A\+{^*}BB\+A^*;$
		\item[(37)]	$(AB)^\dagger=(A^*AB)^\dagger A^*$ and $(A^* AB)^\dagger =B^\dagger(A^*A)^\dagger;$
		\item[(38)]$(AB)^\dagger=B^*(ABB^*)^\dagger$ and $(ABB^*)^\dagger =(BB^*)^\dagger A^\dagger;$
		\item[(39)] $(AB)^\dagger=B^*(A^*ABB^*)^\dagger A^*$ and $(A^*ABB^*)^\dagger =(BB^*)^\dagger(A^* A)^\dagger;$
		\item[(40)] $(AB)^\dagger =(B^*B)^n((AA^*)^mAB(B^*B)^n)^\dagger (AA^*)^m$ and $((AA^*)^mAB(B^*B)^n)^\dagger =(B(B^*B)^n)^\dagger ((AA^*)^mA)^\dagger;$
		\item[(41)] $ (AB)^\dagger=B^*(BB^*)^n ((A^*A)^{m+1}(BB^*)^{n+1})^\dagger (A^*A)^mA^*$ and $((A^*A)^{m+1}(BB^*)^{n+1})^\dagger=((BB^*)^\dagger )^{n+1} ((A^*A)^\dagger)^{m+1}.$
	\end{enumerate}
\end{theorem}	
\begin{proof}\underline{$(1)\Rightarrow(2)$}: Straight forward.\\
	\underline{$(2)\Rightarrow(3)$}: Pre and post-multiplying the given condition by $B^*$ and $A^*,$ respectively, we get $B^*B(AB)^\dagger AA^*=B^*A^*,$ equivalently $AA^*(B^*A^*)^\dagger B^*B=AB.$\\
	\underline{$(3)\Rightarrow(1)$}: We have $B^*B(AB)^\dagger AA^*=B^*A^*.$ Pre- and post-multiplying by $(B^*B)^\dagger$  and $(AA^*)^\dagger $ respectively, we get 
	$B^\dagger B(AB)^\dagger AA^\dagger =B^\dagger A^\dagger.$
It is clear that $ \mathcal R((AB)^\dagger)= \mathcal R((AB)^*) \subseteq R(B^*)$ and $R((AB)^\dagger{^*})=R(AB)\subseteq R(A).$ Thus $B^\dagger B(AB)^\dagger=(AB)^\dagger$ and $ (AB)^\dagger AA^\dagger=(AB)^\dagger.$ Hence, $B^\dagger B(AB)^\dagger AA^\dagger= (AB)^\dagger =B^\dagger A^\dagger.$ \\
	\underline{$(1)\Rightarrow(4)$}: It is easy to see from the assumption that $(AB)^\dagger =(AB)^\dagger AB (AB)^\dagger =B^\dagger A^\dagger ABB^\dagger A^\dagger.$\\
	\underline{$(4)\Rightarrow(5)$}: Pre-multiplying the given condition by $B$ and post-multiplying by 
$A$ we get $B(AB)\+ A=BB\+A\+ABB\+A\+A=(BB\+A\+A)^2.$  $(BB\+A\+A)^4= (BB\+A\+A)^2(BB\+A\+A)^2=B(AB)\+ AB(AB)\+ A=B(AB)\+ A=(BB\+A\+A)^2.$ Since  $BB\+$ and $A\+A$ are orthogonal projections by Lemma \ref{aBBA}, $BB\+A\+A=A\+ABB\+.$ The statements 	$(AB)^\dagger=(A^\dagger AB)^\dagger A^\dagger$ and $(A^\dagger AB)^\dagger =B^\dagger A^\dagger A$ can be proved by verifying all Penrose equations  using  $BB\+A\+A=A\+ABB\+.$\\
	\underline{$(5)\Rightarrow(6)$}: From the assumption, we have $(AB)(AB)^\dagger =AB(A^\dagger AB)^\dagger A^\dagger =ABB^\dagger A^\dagger AA^\dagger =ABB^\dagger A^\dagger.$ Thus, $AB=AB(AB)\+AB=ABB^\dagger A^\dagger AB, (AB)^\dagger=(AB)\+AB(AB)\+=B^\dagger A^\dagger ABB^\dagger A^\dagger,ABB^\dagger A^\dagger $ is a projection and Hermitian. Now $(ABB^\dagger)^\dagger =BB^\dagger A^\dagger $ is easily verifiable, using Theorem \ref{equivalence} $(iii)$. Moreover, $(AB)^\dagger =B^\dagger A^\dagger ABB^\dagger A^\dagger= B^\dagger A^\dagger (AB)(AB)^\dagger =B^\dagger BB^\dagger A^\dagger AB(AB)^\dagger= B^\dagger (ABB^\dagger)^\dagger (AB)(AB)^\dagger=B^\dagger (ABB^\dagger)^\dagger.$\\
	\underline{$(6)\Rightarrow(7)$}:
	Suppose $(AB)^\dagger=B^\dagger(ABB^\dagger)^\dagger.$ Then $(AB)(AB)^\dagger=ABB^\dagger (ABB^\dagger)^\dagger=ABB^\dagger BB^\dagger A^\dagger=ABB^\dagger A^\dagger. $ Thus $AB=ABB^\dagger A^\dagger AB.$ Now by Theorem \ref{equivalence} $(vi),$ we get $(A^\dagger ABB^\dagger)^\dagger =BB^\dagger A^\dagger A.$ 
	 Since $BB^\dagger A^\dagger=(ABB^\dagger)^\dagger ,$ $(A^\dagger ABB^\dagger)^\dagger =(ABB^\dagger)^\dagger A.$ Therefore $B^\dagger(A^\dagger ABB^\dagger)^\dagger A^\dagger =B^\dagger (ABB^\dagger)^\dagger AA^\dagger =(AB)^\dagger AA^\dagger =(AB)^\dagger.$\\
	\underline{$(7)\Rightarrow(8)$}: We have $AB=AB(AB)^\dagger AB=ABB^\dagger (A^\dagger  ABB^\dagger)^\dagger A^\dagger AB=ABB^\dagger BB^\dagger A^\dagger A A^\dagger AB=ABB^\dagger A^\dagger A B$ and $(AB)(AB)^\dagger=ABB^\dagger(A^\dagger ABB^\dagger)^\dagger A^\dagger= ABB^\dagger BB^\dagger A^\dagger AA^\dagger=ABB^\dagger A^\dagger.$ Similarly, $(AB)^\dagger AB=B^\dagger A^\dagger AB.$ Thus $B^\dagger A^\dagger \in AB\{1,3,4\}.$\\
		\underline{$(8)\Rightarrow(9)$}: Since $B^\dagger A^\dagger \in AB\{1,3,4\},$  $ AB=ABB\+A\+AB$ and $ ABB\+A\+=(ABB\+A\+)^*.$ Now, $(AB)(AB)\+=ABB\+A\+AB(AB)\+=(ABB\+A\+)^*AB(AB)\+=(B\+A\+)^*(AB)^*AB(AB)\+=(B\+A\+)^*(AB)^*=(ABB\+A\+)^*=ABB\+
		A\+=A\+{^*}BB\+A^*.$ Similarly we can prove the  other relation.\\
		\underline{$(9)\Rightarrow(10)$}: Since $AB(AB)\+ =ABB\+A\+,$  $AB(AB)\+ AB=ABB\+ A\+ AB.$ Then by Theorem \ref{equivalence} $(iv)$, $A\+ ABB\+ =BB\+A\+ A.$ It is clear from the assumption, we have $B\+ A^*\in A{\+}^*B\{3,4\}.$
 Also, it is easy to verify that $B\+A^*\in A{\+}^*B\{1,2\}.$	\\	
		\underline{$(10)\Rightarrow(1)$}: Applying Theorem \ref{equivalence} for $A\+{^*}$ and $B,$ we get $A^*A\+{^*}BB\+=BB\+A^*A\+{^*} $ i.e., $A\+ABB\+=BB\+A\+A.$ Using the third and fourth Penrose conditions for (10), we get
\begin{align*}
     ABB\+A\+AB&=  AA\+ABB\+B=AB,\\
 B\+A\+ABB\+A\+&=B\+BB\+A\+AA\+=B\+A\+,\\(ABB\+A\+)^*&=A\+{^*}BB\+A^*=ABB\+A\+,\\(B\+A\+AB)^*&=B^*A\+AB\+{^*}=B\+A\+AB.
\end{align*}
The equivalences of (10)-(18) can be established by replacing $A$ by $A\+{^*}$ in (1)-(9).
Similarly, the equivalences of (19)-(27) can be established by replacing $B$ by $B\+{^*}$ in (1)-(9) and
the equivalences of (28)-(36) can be established by replacing $A$ by $B\+$ and $B$ by $A\+$ in (1)-(9). The equivalence of (1) and (19) is similar to that of (1) and (10). The equivalence of (1) and (28) follows by applying Moore-Penrose inverse on both sides of (1) and (28).  \\
	\underline{$(1)\Rightarrow(37)$}: We use Theorem \ref{equivalence} to verify Penrose equations to prove $(A^* AB)^\dagger =B^\dagger(A^*A)^\dagger.$ We get the first Penrose equation 
	verified as below. 
\begin{align*} A^* AB B^\dagger(A^*A)^\dagger A^* AB= A^* AB B^\dagger A^\dagger  AB=  A^* AA^\dagger AB B^\dagger B=A^* AB .\end{align*} By Theorem \ref{equivalence}  and Theorem \ref{equivalence123} $(v)$, we get
	\begin{align*}
	   [A^* AB B^\dagger(A^*A)^\dagger]^*&=[BB^\dagger A^*A(A^*A)^\dagger]^*= (BB^\dagger A^\dagger A)^*=A^\dagger ABB^\dagger .
	\end{align*}The right side is Hermitian, so is the left side.
	Similarly, we can prove the second and fourth Penrose equations.
Also we get	\begin{align*}(A^*AB)^\dagger A^* =B^\dagger(A^*A)^\dagger A^* =B\+A\+=(AB)^\dagger\end{align*}  \\
	\underline{$(1)\Rightarrow(38)$}: Similar to $(1)\Rightarrow(37).$\\
	\underline{$(1)\Rightarrow(39)$}: By using Theorem \ref{equivalence}, we get \begin{align*}A^*ABB^* (BB^*)^\dagger(A^* A)^\dagger A^*ABB^*&=A^*ABB^\dagger A^\dagger ABB^*\\&=A^*A A^\dagger ABB^\dagger BB^*=A^*ABB^*. \end{align*} Now using  Theorem \ref{equivalence123} $(v)$,
\begin{align*}  [A^*ABB^* (BB^*)^\dagger(A^* A)^\dagger]^*&=[A^*ABB^\dagger(A^* A)^\dagger]^*\\& =[BB^\dagger A^*A(A^* A)^\dagger]^*\\&=(BB^\dagger A^\dagger A)^*,
\end{align*} which is Hermitian by Theorem \ref{equivalence}. Hence the first and  third conditions of Penrose equations are satisfied. The second and fourth conditions  follow similarly.\\
\underline{$(1)\Rightarrow(40)$}: Let $P=(AA^*)^mA$ and $Q=B(B^*B)^n. $ Then they have closed ranges since by Lemma \ref{closedrangelemma}, $\mathcal R((AA^*)^m A)$  is closed if and only if $\mathcal R((AA^*)^m{\+} (AA^*)^m AA\+)$ is closed. But the latter is closed as it equals $\mathcal R (AA\+)=\mathcal R (AA^*),$ whose closedness follows from the closedness of $ \mathcal R(A).$ Similar argument works for $Q$ also. Now again by Lemma \ref{closedrangelemma}, $PQ$ has a closed range if and only if $P\+PQQ\+$ has a closed range. For,
 \begin{align*}P\+PQQ\+&=[(AA^*)^mA]\+(AA^*)^mA B(B^*B)^n[B(B^*B)^n]\+\\&=A\+[(AA^*){^m} ]\+(AA^*)^m A B(B^*B)^n[(B^*B){^n}]\+ B\+\\&=A\+AA\+ABB\+BB\+ =A\+ABB\+.
 \end{align*} Here the reverse order law is applied for $[(AA^*)^mA]\+$ and $[B(B^*B)^n]\+$
 as they satisfy the condition $(ii)$  in Theorem \ref{Rolrange}. This ensures the existence of  $((AA^*)^mAB(B^*B)^n)\+.$

We prove $[(AA^*)^mAB(B^*B)^n]^\dagger =[B(B^*B)^n]^\dagger [(AA^*)^mA]^\dagger$ i.e., $(PQ)\+=Q\+P\+$ by verifying the Penrose equations. We have \begin{align*}
PQQ\+ P\+ PQ&=(AA^*)^m AB(B^*B)^n [(B^*B)^n] {\+} B\+ A\+ [(AA^*)^m] {\+} (AA^*)^m AB(B^*B)^n\\&=(AA^*)^m ABB\+  BB\+A\+ AA\+ AB(B^*B)^n\\&=(AA^*)^m ABB\+ A\+ AB(B^*B)^n \\&=(AA^*)^m AA\+ ABB\+ B(B^*B)^n \\&=(AA^*)^m AB(B^*B)^n=PQ.
\end{align*}
Similarly, we can prove the second Penrose equation. For proving the third one we use the following facts $A\+ (AA^*)^m=A^*(AA^*)^{m-1},$ $[(AA^*)^m]\+ A=[(AA^*)^{m-1}]\+ A\+{^*}$ for all $m\ge1$ and $(ABB\+ A\+)^*=ABB\+ A\+.$ We have
\begin{align*}
(PQQ\+P\+)^*&=((AA^*)^m ABB\+ A\+ [(AA^*)^m]\+)^*\\&= [(AA^*)^m]\+ABB\+ A\+(AA^*)^m\\&=A^*{\+}BB\+ A^*\\&=(ABB\+A\+)^*.
\end{align*}
The right side is Hermitian so is the left side. Similarly, the fourth Penrose equation can be proved. Also, we have
\begin{align*}
(B^*B)^n[(AA^*)^mAB(B^*B)^n]^\dagger (AA^*)^m&=(B^*B)^n[B(B^*B)^n]^\dagger [(AA^*)^mA]\+(AA^*)^m\\&=(B^*B)^n[(B^*B)^n]{\+}B\+ A\+[(AA^*)^m]{\+}(AA^*)^m\\&=B\+BB\+A\+AA\+=B\+A\+=(AB)\+.
\end{align*}\\
	\underline{$(1)\Rightarrow(41)$}: Let $P=(AA^*)^{m+1}$ and $Q=(B^*B)^{n+1}.$ We can prove the existence of $(PQ)\+$ using Lemma \ref{closedrangelemma}, by a similar argument in $(1)\Rightarrow(40)$.
	\begin{align*}
	PQQ\+P\+PQ&=(AA^*)^{m+1}B\+BAA\+(B^*B)^{n+1}\\&=(AA^*)^{m+1}AA\+B\+B(B^*B)^{n+1}\\&=(AA^*)^{m+1}(B^*B)^{n+1}=PQ.
	\end{align*}
	Using the fact $AA^*B\+B=B\+BAA^*,$ we have
	\begin{align*}
	(PQQ\+P\+)^*&=(AA^*)^{m+1}B\+B[(AA^*)^{m+1}]{\+}\\&=B\+B(AA^*)^{m+1}[(AA^*)^{m+1}]{\+}\\&=B\+BAA\+.
		\end{align*}
	The right side is Hermitian so is the left side. Similarly, the second and  fourth Penrose equations can be proved. Now we have
	\begin{align*}
B^*(BB^*)^n ((A^*A)^{m+1}(BB^*)^{n+1})^\dagger (A^*A)^mA^*&=B^*(BB^*)^n[(BB^*)^{n+1}]{\+}[(A^*A)^{m+1}]{\+}\\  &\hspace{5cm}(A^*A)^mA^*\\&=B^*BB\+(BB^*)\+(A^*A)\+A\+AA^*\\&=B^*(BB^*)\+(A^*A)\+A^*\\&=B\+A\+=(AB)\+.
		\end{align*}\\
	\underline{$(37)-(41)\Rightarrow(1)$}: Using the given conditions we get	\begin{align*}(AB)\+=(A^*AB)^\dagger A^* =B^\dagger(A^*A)^\dagger A^*=B\+A\+.
	\end{align*} 
	Similarly, other implications also follow by substituting the second set of equations into the first ones and using the properties of Moore-Penrose inverse  given in Lemma \ref{properties}.
\end{proof}

\begin{center}
	{\bf Acknowledgement}
\end{center}
The first author thanks National Institute of Technology Karnataka (NITK), Surathkal for the financial support.

\end{document}